\documentclass{amsart}

\allowdisplaybreaks
\usepackage[all]{xy}

\newcommand{\F}{\mathcal{F}}

\newcommand{\G}{\mathcal{G}}
\newcommand{\Pc}{\mathcal{P}}
\newcommand{\Hc}{\mathcal{H}}

\newcommand{\bN}{\mathbb{N}}
\newcommand{\bP}{\mathbb{P}}
\newcommand{\bQ}{\mathbb{Q}}
\newcommand{\bR}{\mathbb{R}}
\newcommand{\bT}{\mathbb{T}}
\newcommand{\bZ}{\mathbb{Z}}

\DeclareMathOperator{\tr}{tr}

\theoremstyle{plain} \newtheorem{theorem}{Theorem}[section]
\theoremstyle{plain} \newtheorem{proposition}[theorem]{Proposition}
\theoremstyle{plain} \newtheorem{lemma}[theorem]{Lemma}
\theoremstyle{plain} \newtheorem{corollary}[theorem]{Corollary}
\theoremstyle{definition} \newtheorem{definition}[theorem]{Definition}
\theoremstyle{definition} \newtheorem{notation}[theorem]{Notation}
\theoremstyle{remark} \newtheorem{remark}[theorem]{Remark}

\numberwithin{equation}{section}

\begin{document}

\title[stochastic equations]{Stochastic equations \\
on projective systems of groups}

\author{Steven N.\ Evans}
\address{Steven N.\ Evans \\
  Department of Statistics \#3860 \\
  University of California at Berkeley \\
  367 Evans Hall \\
  Berkeley, CA 94720-3860 \\
  U.S.A.}
\email{evans@stat.Berkeley.EDU}
\urladdr{http://www.stat.berkeley.edu/users/evans/}

\author{Tatyana Gordeeva}
\address{ Tatyana Gordeeva \\
  Department of Statistics \#3860 \\
  University of California at Berkeley \\
  367 Evans Hall \\
  Berkeley, CA 94720-3860 \\
  U.S.A.}
\email{gordeeva@stat.Berkeley.EDU}
\urladdr{http://www.stat.berkeley.edu/users/gordeeva/}

\thanks{SNE supported in part by NSF grant DMS-0907630. TG supported in part
by a VIGRE grant awarded to the Department of Statistics, University of California
at Berkeley}


\begin{abstract}
We consider stochastic equations of the form $X_k = \phi_k(X_{k+1}) Z_k$, $k \in \mathbb{N}$,
where $X_k$ and $Z_k$ are random variables taking values in a compact group $G_k$,
$\phi_k: G_{k+1} \to G_k$ is a continuous homomorphism, and the noise
$(Z_k)_{k \in \mathbb{N}}$ is a sequence of independent random variables. We take
the sequence of homomorphisms and the sequence of noise distributions as given,
and investigate what conditions on these objects result in a unique distribution
for the ``solution'' sequence $(X_k)_{k \in \mathbb{N}}$ and what conditions
permits the existence of a solution sequence that is a function of the noise alone 
(that is, the solution does not incorporate extra input randomness ``at infinity'').
Our results extend previous work on stochastic equations on a single group
that was originally motivated by Tsirelson's example of a stochastic differential 
equation that has a unique solution in law but no strong solutions. 
\end{abstract}

\maketitle

\section{Introduction}

The following stochastic process was considered by Yor in \cite{MR1147613}
in order to clarify the structure underpinning
Tsirelson's celebrated example \cite{MR0375461} of a stochastic differential equation 
that does not have a strong solution even though all solutions have
the same law.  

Let $\bT$ be the usual circle group; that is, $\bT$
can be thought of as the interval $[0,1)$ equipped with addition modulo
$1$.    Suppose for each $k \in \bN$ that 
$\mu_k$ is a Borel probability measure on $\bT$. 
Write $\mu = (\mu_k)_{k \in \bN}$.
We say that sequence of $\bT$-valued random variables $(X_k)_{k \in \bN}$
defined on some probability space $(\Omega, \F, \bP)$ 
{\em solves the stochastic equation associated with}
$\mu$ if 
\[
\bP[f(X_k) \, | \, (X_j)_{j>k}]
= \int_\bT f(X_{k+1} + z) \, \mu_k(dz)
\]
for all bounded Borel function $f: \bT \to \bR$, where
we use the notation $\bP[ \cdot \, | \, \cdot ]$ for condition expectations
with respect to $\bP$.  In other words, if for each
$k \in \bN$ we define
a $\bT$-valued random variable $Z_k$ by requiring
\begin{equation}
\label{Yor_example}
X_k = X_{k+1} + Z_k,
\end{equation}
then $(X_k)_{k \in \bN}$ solves the stochastic
equation associated with $\mu$ if and only if for all $k \in \bN$
the distribution of $Z_k$ is $\mu_k$ and $Z_k$ is independent of
$(X_j)_{j > k}$.

Yor addressed the existence of solutions $(X_k)_{k \in \bN}$
that are {\em strong} in the sense that the random variable
$X_k$ is measurable with respect to
$\sigma((Z_j)_{j \ge k})$ for each $k \in \bN$;
that is, speaking somewhat informally,
a solution is strong if it can be reconstructed from the
``noise'' $(Z_j)_{j \in \bN}$ without introducing
additional randomness ``at infinity.'' It turns out that
strong solutions exist if and only if
\[
\lim_{m \to \infty} \lim_{n \to \infty} 
\prod_{\ell=m}^n
\left|\int_\bT \exp(2 \pi i h x) \, \mu_\ell(dx)\right|
> 0
\]
for all $h \in \bZ$ or, equivalently,
\[
\sum_{k=1}^\infty 
\left[1 - \left|\int_\bT \exp(2 \pi i h x) \, \mu_k(dx)\right| \right]
< \infty.
\]

Yor's investigation was extended in \cite{MR2365485}, where the group
$\bT$ is replaced by an arbitrary, possibly non-abelian, compact
Hausdorff group.  As one would expect, the role of the the complex exponentials
$\exp(2 \pi i h \cdot)$, $h \in \bZ$, in this more
general setting is played by group representations.
Interesting new phenomena appear when the group is non-abelian
due to the fact that there are irreducible representations which
are no longer one-dimensional.  Several of the
 results in \cite{MR2365485} are framed in terms of properties
 of the set of extremal solutions (that is, solutions that
 can't be written as mixtures of others), and the structure of
 such solutions was elucidated further in \cite{MR2653259}.

We further extend the work in \cite{MR1147613, MR2365485} by considering the
following more general set-up.

Fix a sequence $(G_k)_{k \in \bN}$ 
of compact Hausdorff groups with countable bases. Suppose for each $k \in \bN$
that there is a continuous homomorphism $\phi_k: G_{k+1} \to G_k$.
Define a compact subgroup $H \subseteq G := \prod_{k \in \bN} G_k$
by
\begin{equation}
\label{defH}
H:= \{g = (g_k)_{k \in \bN} \in G : g_k = \phi_k(g_{k+1}) \, \text{for all $k \in \bN$}\},
\end{equation}

For example, if we take $G_k = \bT$ for all $k \in \bN$, then
the homomorphism $\phi_k$ is necessarily of the form $\phi_k(x) = N_k x$
for some $N_k \in \bZ$ and 
\[
H = \{g = (g_k)_{k \in \bN} \in G : g_k = N_k g_{k+1} \, \text{for all $k \in \bN$}\}.
\]

For a more interesting example, fix a compact group abelian group $\Gamma$,
put  $G_k := G_{1,k} \times G_{2,k-1} \cdots \times G_{k,1}$, where
each group $G_{i,j}$ is a copy of $\Gamma$, and define the homomorphism
$\phi_k$ by
\[
\phi_k(g_{1,k+1}, g_{2,k}, \ldots, g_{k+1,1}) 
:= (g_{1,k+1} + g_{2,k}, g_{2,k}+g_{3,k-1}, \ldots, g_{k,2} + g_{k+1,1})
\]
(where we write the group operation in $\Gamma$ additively).
Note that in this case $H$ is isomorphic to the infinite product $\Gamma^{\bN}$, 
because an element $h = (h_{i,j})_{(i,j) \in \bN \times \bN}$
is uniquely specified by the values  $(h_{i,1})_{i \in \bN}$
and there are no constraints on these elements.
The following pictures shows a piece of an element of
$H$ when $\Gamma$ is the group $\{0,1\}$ equipped with addition modulo $2$.

\centerline{
\xymatrix @R=2pc @C=2pc{
*+[Fo]{1}\ar[d]\\
*+[Fo]{0}\ar[d] & *+[Fo]{1}\ar[l]\ar[d]\\
*+[Fo]{1}\ar[d] & *+[Fo]{1}\ar[l]\ar[d] & *+[Fo]{0}\ar[l]\ar[d]\\
*+[Fo]{1}\ar[d] & *+[Fo]{0}\ar[l]\ar[d] & *+[Fo]{1}\ar[l]\ar[d] & *+[Fo]{1}\ar[l]\ar[d]\\
*+[Fo]{1} & *+[Fo]{0}\ar[l] & *+[Fo]{0}\ar[l] & *+[Fo]{1}\ar[l] & *+[Fo]{0}\ar[l]\\
}
}

Assume for each $k \in \bN$ that $\mu_k$ is a Borel probability measure 
$G_k$ and write $\mu = (\mu_k)_{k \in \bN}$.
We say that sequence of random variables $(X_k)_{k \in \bN}$
defined on some probability space $(\Omega, \F, \bP)$, where
$X_k$ takes values in $G_k$, {\em solves the stochastic equation associated with}
$\mu$ if 
\[
\bP[f(X_k) \, | \, (X_j)_{j>k}]
= \int_{G_k} f(\phi_k(X_{k+1}) z) \, \mu_k(dz)
\]
for all bounded Borel function $f: G_k \to \bR$.  In other words, if for each
$k \in \bN$ we define
a $G_k$-valued random variable $Z_k$ by requiring
\begin{equation}\label{sde}
X_k = \phi_k(X_{k+1}) Z_k,
\end{equation}
then $(X_k)_{k \in \bN}$ solves the stochastic
equation if and only if for all $k \in \bN$
the distribution of $Z_k$ is $\mu_k$ and $Z_k$ is independent of
$(X_j)_{j > k}$.  In particular, if $(X_k)_{k \in \bN}$ solves the stochastic
equation, then the sequence of random variables
$(Z_k)_{k \in \bN}$ is independent. 

Certain special cases of this set-up when $G_k = \Gamma$, $k \in \bN$, 
for some fixed group $\Gamma$ and  $\phi_k = \psi$, $k \in \bN$ 
for a fixed automorphism $\psi: \Gamma \to \Gamma$  were considered in
\cite{MR2582432, Raja}.

Note that whether or not a sequence 
$(X_k)_{k \in \bN}$ solves the stochastic equation associated with 
$\mu$ is solely a feature 
of the distribution of the sequence, and so we say that a probability
measure on the product group
$\prod_{k \in \bN} G_k$ is a solution of the stochastic equation
if it is the distribution of a sequence that solves the equation
and write $\Pc_\mu$ for the set of such measures.

In keeping with the terminology above, we say that a solution  $(X_k)_{k \in \bN}$
is {\em strong} if $X_k$ is measurable with respect to
$\sigma((Z_j)_{j \ge k})$ for each $k \in \bN$.
Note that whether or not a solution is strong also depends only its distribution,
and so we define strong elements of $\Pc_\mu$ in the obvious manner
and denote the set of such probability measures by $\Pc_\mu^{\mathrm{strong}}$.

Because
applying the homomorphism 
$\phi_k$ to $X_{k+1}$  can degrade the ``signal''
present in $X_{k+1}$ (for example, $\phi_k$ need not
be invertible), 
the question of whether or not strong solutions exist 
will involve the interaction between the homomorphisms $(\phi_k)_{k \in \bN}$
and distributions $(\mu_k)_{k \in \bN}$ of the noise random variables
and it introduces new phenomena not present in \cite{MR1147613, MR2365485}.

An outline of the rest of the paper is as follows.
In the Section \ref{S:extreme_and_strong} we examine the compact, convex set of solutions
and show that strong solutions are extreme points of this set. 
We show that the subgroup $H$ acts transitively on the extreme points of the set of solutions
and we relate the existence of strong solutions to properties
of the set of extreme points. In Section \ref{S:representations_and_strong},
 we obtain criteria for the existence of strong solutions in terms of the 
the representations of the group $G_k$ and the corresponding
Fourier transforms of the probability measures $\mu_k$. 
In Section \ref{S:representations_and_strong}, we determine the relationship
between the existence of strong solutions and the phenomenon of ``freezing''
wherein almost all sample paths of the  random noise sequence 
agrees with some sequence of
constants for all sufficiently large indices. 
Finally, in Section \ref{S:lattice} and \ref{S:torus}, respectively,
we investigate the example considered above of
random variables indexed by the nonnegative quadrant of the two-dimensional
integer lattice and another example where each group $G_k$ is
the two dimensional torus and each homomorphisms $\phi_k$
is a fixed ergodic toral automorphism.

\section{Extreme points of $\Pc_\mu$ and strong solutions}
\label{S:extreme_and_strong}

It is natural to first inquire whether $\Pc_\mu$ is non-empty and, if so,
whether it consists of a single point; that is, whether there exist
probability measures that solve the stochastic equation associated with
$\mu$ and, if so, whether there is a single such measure.  The question
of existence is easily disposed of by Proposition~\ref{existsolution}
below.  Note that
because the group $G = \prod_{k \in \bN} G_k$ is compact and metrizable,
the set of probability measures on $G$ equipped with the topology of 
weak convergence is also compact and metrizable.

\begin{proposition} 
\label{existsolution}
For any sequence $\mu$, the set $\Pc_\mu$ is non-empty.
\end{proposition}
\begin{proof}
Construct on some probability space  
a sequence $(Z_k)_{k \in \bN}$ of independent  random variables
such that $Z_k$ has distribution $\mu_k$.
For each $N \in \bN$, define random variables 
$X^{(N)}_1, \ldots, X^{(N)}_{N+1}$ recursively by
\[
X^{(N)}_{N+1} := e_{N+1} := \text{identity in $G_{N+1}$}
\]
and
\[
X^{(N)}_{k} = \phi_{k}(X^{(N)}_{k+1}) Z_k, \quad 1 \leq k \leq N,
\]
so that for $1 \le k \le N$ the
random variable $\phi_k(X_{k+1}^{(N)})^{-1} X_k^{(N)} $ has distribution $\mu_k$ 
and is independent of $X_{k+1}^{(N)}, X_{k+2}^{(N)}, \ldots, X_{N}^{(N)}$.

Write $\bP_N$ for the distribution of 
the sequence $(X_1^{(N)}, \ldots, X_N^{(N)}, e_{N+1}, e_{N+2}, \ldots)$. 
Because the space of probability measures
on the group $\prod_{k \in \bN} G_k$ equipped with the
weak topology is compact and metrizable, there 
exists a subsequence $(N_n)_{n \in \bN}$ and a probability measure $\bP_\infty$ 
such that $\bP_{N_n} \to \bP_\infty$ weakly as $n \to \infty$. It is clear that 
$\bP_\infty \in \Pc_\mu$.
\end{proof}

The question of uniqueness (that is, whether or not $\# \Pc_\mu = 1$)
 is more demanding and will occupy much of our 
attention in the remainder of the paper.

As a first indication of what is involved, consider the case where
each measure $\mu_k$ is simply the unit point mass at the identity $e_k$ of $G_k$.
In this case $(X_k)_{k \in \bN}$ solves the stochastic equation if 
$X_k = \phi_k(X_{k+1})$ for all $k \in \bN$.  
Recall the definition of the compact subgroup 
$H \subseteq G := \prod_{k \in \bN} G_k$
from \eqref{defH}.
It is clear that $\Pc_\mu$ coincides with the set of probability measures that
are supported on $H$, and hence $\# \Pc_\mu = 1$  if and only
if $H$ consists of just the single identity element.  Note that
if $\# H > 1$  and $(X_k)_{k \in \bN}$ is a solution
with distribution $\bP \in \Pc_\mu$ that is not a point mass, then
$X_k$ is certainly not a function of $(Z_j)_{j \ge k}
= (e_j)_{j \ge k}$ and the solution $(X_k)_{k \in \bN}$
is not strong.  Moreover, the probability measures $\bP \in \Pc_\mu$
that are distributions of strong solutions $(X_k)_{k \in \bN}$
are the point masses at elements of $H$ and $\Pc_\mu$ is the closed convex hull
of this set of measures.

An elaboration of the argument we have just given establishes the following
result.

\begin{proposition}
\label{existsnonstrongsolution} 
If $H$ is non-trivial (that is, contains elements other than the identity),
then $\Pc_\mu \setminus \Pc_\mu^{\mathrm{strong}} \ne \emptyset$.
In particular, if $H$ is non-trivial and $\# \Pc_\mu = 1$, then
$\Pc_\mu^{\mathrm{strong}} = \emptyset$.
\end{proposition}

\begin{proof} Suppose that all solutions are strong.
Let $(X_k)_{k \in \bN}$ be a strong solution.
 
By extending the underlying probability space if necessary,
construct an $H$-valued random variable $(U_k)_{k \in \bN}$ that is independent
of $(X_k)_{k \in \bN}$ and is not almost surely constant. 
Note that $(U_k)_{k \in \bN}$ is not $\sigma((X_k)_{k \in \bN})$-measurable and hence,
{\em a fortiori}, $(U_k)_{k \in \bN}$ is not $\sigma((Z_k)_{k \in \bN})$-measurable.

Observe that
\[\phi_k(U_{k+1}X_{k+1}) Z_k = \phi_k(U_{k+1})\, \phi_k(X_{k+1}) Z_k = U_k X_{k},\]
because $\phi_k(U_{k+1}) = U_k$ for all $k \in \bN$ by definition of $H$.
Hence, $(U_k X_k)_{k \in \bN}$ is also a solution. 
Thus,  $(U_k X_k)_{k \in \bN}$ is a strong solution by our assumption that 
all solutions are strong. In particular,
$U_k X_k$ is  $\sigma((Z_j){j \geq k})$-measurable for all $k \in \bN$.
However, $U_k = (U_k X_k) X_k^{-1}$ is $\sigma((Z_j)_{j \geq k})$-measurable, 
and we arrive at a contradiction.
\end{proof}

\begin{remark}
Consider  the particular setting of \cite{MR2365485}, 
where $G_k = \Gamma$, $k \in \bN$, for
some fixed group $\Gamma$, each homomorphism $\phi_k$ is the identity,
and $H = \{(g,g,\ldots) : g \in \Gamma\}$.  In this case, 
one can choose the sequence $(U_k)_{k \in \bN}$ in the proof of
Proposition~\ref{existsnonstrongsolution} to be $(U,U,\cdots)$,
where $U$ is distributed according to Haar measure on $\Gamma$; that is,
$(U_k)_{k \in \bN}$ is distributed according to Haar measure on $H$.
Each marginal distribution of the solution $(X_k)_{k \in \bN}$ 
is then Haar measure on $G_k = \Gamma$.  In our more general setting
it will not generally be the case that if $(U_k)_{k \in \bN}$ is
distributed according to Haar measure on $H$, then $X_k$ will
be distributed according to Haar measure on $G_k$ for each $k \in \bN$.
For example, fix a compact group $\Gamma$,
put $G_k = \Gamma^{\bN}$ for all $k \in \bN$ and define
$\phi_k : G_{k+1} \to G_k$ by $\phi_k(g_1,g_2,g_3, \ldots)
= (g_1,g_1,g_2,g_2, g_3, g_3, \ldots)$ for all $k \in \bN$.
It is clear that 
$H = \{((g,g,\ldots), (g,g,\ldots), \ldots) : g \in \Gamma\}$, so that
$\{x_k : (x_1, x_2, \ldots) \in H\} \subseteq G_k$
is just the diagonal subgroup
$\{(g,g,\ldots) : g \in \Gamma\}$ of the group $G_k$.  Hence, for example,
if $\mu_k$ is the point mass at the identity of $G_k$ for each $k \in \bN$,
the possible solutions $(X_k)_{k \in \bN}$ are just
arbitrary random elements of $H$, and it is certainly not possible
to construct a solution such that the marginal distribution of $X_k$
is Haar measure on $G_k$ for some $k \in \bN$.
\end{remark}

From now on, we  let $X_k : G \to G_k$, $k \in \bN$, denote the random variable
defined by $X_k((x_j)_{j \in \bN}) := x_k$ and define $Z_k : G \to G_k$, $k \in \bN$,
by $Z_k := \phi_n(X_{k+1})^{-1} X_k$.

\begin{notation} Given a sequence of random variables 
$S = (S_1, S_2, \ldots)$ and $k \in \bN$, set $\F^S_k := \sigma((S_j)_{j \ge k})$.
Similarly, set $\F^S := \F_1^S$ and $\F^S_\infty := \bigcap_{k \in \bN} \F^S_k$.
\end{notation}

\begin{notation} For any sequence $\mu = (\mu_k)_{k \in \bN}$,
the set of solutions $\Pc_\mu$ is clearly a compact convex subset. 
Let $\Pc_\mu^{\mathrm{ex}}$ denote the extreme points of $\Pc_\mu$.
\end{notation}

\begin{lemma}
\label{extremetrivial} A probability measure $\bP \in \Pc_\mu$
belongs to $\Pc_\mu^{\mathrm{ex}}$if and only if the remote future 
$\F_{\infty}^{X}$ is trivial under $\bP$.
\end{lemma}

\begin{proof}
Our proof follows that of an analogous result in \cite{MR2365485}. 

Suppose that $\bP \in \Pc_\mu$ and the $\sigma$-field $\F_\infty^{X}$ is not trivial under $\bP$. 

Fix a set $A \in \F_\infty^X$ with $0 < \bP(A) < 1$. Then, 
\[
\bP(\cdot) = \bP(A)\bP(\cdot \, | \, A) + (1 - \bP(A))\bP(\cdot \, | \, A^c).
\]
Observe that $\bP(\cdot \, | \, A) \neq \bP(\cdot \, | \, A^c)$, 
since $\bP(A \, | \, A) = 1 \neq \bP(A \, | \, A^c) = 0$. 

Note for each $k \in \bN$ and $B \subseteq G_k$ that
\begin{align*}
\bP\{X_k \, \phi_k(X_{k+1})^{-1} \in B \, | \, A\} 
&= \frac{\bP(\{X_k \, \phi_k(X_{k+1})^{-1} \in B\} \cap A)}{\bP(A)}\\
& = \frac{\mu_k(B)\bP(A)}{\bP(A)} = \mu_k(B)
\end{align*}
 because $\bP \in \Pc_\mu$ and hence  $X_k \, \phi_k(X_{k+1})^{-1}$ 
 is independent of $\F_\infty^{X}$ under $\bP$. 
Similarly, if $C \in \F^{X}_{k+1}$, 
\begin{align*}
\bP(\{X_k \, \phi_k(X_{k+1})^{-1} \in B\} \cap C \, | \, A) 
&= \frac{\mu_k(B)\bP(C \cap A)}{\bP(A)}\\ 
&= \bP\{X_k \, \phi_k(X_{k+1})^{-1} \in B \, | \, A\}\bP(C \, | \, A)
\end{align*}

Thus, $\bP(\cdot \, | \, A) \in \Pc_\mu$. 
The analogous argument establishes
$\bP(\cdot \, | \, A^c) \in \Pc_\mu$.
Since $\bP(\cdot \, | \, A) \ne \bP(\cdot \, | \, A^c)$,
the probability measure $\bP$ cannot belong to $\Pc_\mu^{\mathrm{ex}}$.

Now assume that 
$\bP \in \Pc_\mu$ and $\F_\infty^{X}$ is trivial under $\bP$. 
To show $\bP$ is an extreme point, it suffices to show that if $\bP' \in \Pc_\mu$ 
is absolutely continuous with respect to $\bP$, then $\bP = \bP'$. 

Note that a solution $X$ is a time-inhomogeneous Markov chain 
(indexed in backwards time with index set starting at infinity) 
with the following transition probability:
\[
\bP\{X_k \in A \, | \, X_{k+1}\} = \mu_k\{g \in G_k : \phi_k(X_{k+1})g \in A\}.
\]

Since $\bP$ and $\bP'$ are 
the distributions of Markov chains with common transition probabilities 
and $\bP'$ is absolutely continuous with respect to $\bP$, 
it follows that
for any measurable set $A$ the random variables 
$\bP(A \, | \, \F_\infty^{X})$ and $\bP'(A  \, | \, \F_\infty^{X})$ are equal 
$\bP$-a.s. Because $\F_\infty^{X}$ is trivial under both $\bP$ and $\bP'$, 
it must be the case that $\bP(A) = \bP'(A)$.
\end{proof}

\begin{corollary}
\label{strong_are_extreme}
All strong solutions $\bP \in \Pc_\mu$ are extreme; that is,
$\Pc_\mu^{\mathrm{strong}} \subseteq \Pc_\mu^{\mathrm{ex}}$.
\end{corollary}

\begin{proof} 
By definition, if $\bP \in \Pc_\mu$ is strong, then
$X_k \in \F_k^Z$ for all $k \in \bN$.  Thus, $\F_k^X = \F_k^Z$ for all $k \in \bN$
and hence $\F_\infty^X = \F_\infty^Z$.  The last $\sigma$-field 
is trivial by the Kolmogorov zero-one law.
\end{proof}

\begin{remark}
There can be extreme solutions that are not strong.  For example, suppose
that the $G_k = \Gamma$, $k \in \bN$, for some non-trivial group $\Gamma$,
each $\phi_k$ is the identity map, and each $\mu_k$ is the Haar measure on $\Gamma$.
It is clear that $\Pc_\mu$ consists of just the measure 
$\bigotimes_{k \in \bN} \mu_k$ (that is, Haar measure on $G$), and so
this solution is extreme.  However, it follows from
Proposition~\ref{existsnonstrongsolution} that this solution is not strong.
\end{remark}

It is clear that if $\bP \in \Pc_\mu$ and $h = (h_k)_{k \in \bN} \in H$, then 
the distribution of the sequence $(h_k X_k)_{k \in \bN}$ also belongs to 
$\bP \in \Pc_\mu$.
Moreover, if $\bP \in \Pc_\mu^{\mathrm{ex}}$, then it follows from 
Lemma~\ref{extremetrivial} that the distribution of the sequence
$(h_k X_k)_{k \in \bN}$ 
also belongs to $\Pc_\mu^{\mathrm{ex}}$.  
Similarly, if $\bP \in \Pc_\mu^{\mathrm{strong}}$, 
then the distribution of the sequence $(h_k X_k)_{k \in \bN}$
also belongs to $\Pc_\mu^{\mathrm{strong}}$.  We record these observations for 
future reference.

\begin{lemma} The collection of maps $T_h: \Pc_\mu \to \Pc_\mu$, $h \in H$, 
defined by
 $T_h(\bP)(\cdot) = \bP\{(h_k X_k)_{k \in \bN} \in \cdot\}$ constitute a 
a group action of $H$ on $\Pc_\mu$.   The set $\Pc_\mu^{\mathrm{ex}}$ of 
extreme solutions
and the set $\Pc_\mu^{\mathrm{strong}}$ of strong solutions 
are both invariant for this action.
\end{lemma}

It follows from the next result that either 
$\Pc_\mu^{\mathrm{strong}} = \emptyset$
or $\Pc_\mu^{\mathrm{strong}} = \Pc_\mu^{\mathrm{ex}}$.  For the purposes
of the proof and later it is  convenient to introduce the following notation.

\begin{notation}
For $k,\ell \in \bN$ with $k < \ell$, define $\phi^\ell_k : G_\ell \to G_k$ by
\[
\phi_k^\ell= \phi_{k} \circ \phi_{k+1} \circ \cdots \circ \phi_{\ell-1},
\]
and adopt the convention that $\phi_k^k$ is the identity map from $G_k$
to itself.
\end{notation}

\begin{theorem}
\label{htransitive} 
The group action $(T_h)_{h \in H}$ is transitive on $\Pc_\mu^{\mathrm{ex}}$.
\end{theorem}

\begin{proof}
For $k \in \bN$, define 
$X_k' : \prod_{k \in \bN} (G_k \times G_k \times G_k) \to G_k$
(resp. $X_k'' : \prod_{k \in \bN} (G_k \times G_k \times G_k) \to G_k$)
and $Y_k : \prod_{k \in \bN} (G_k \times G_k \times G_k) \to G_k$)
by $X_k'((x_j',x_j'',y_j)_{j \in \bN}) = x_k'$ (resp. 
$X_k''((x_j',x_j'',y_j)_{j \in \bN}) = x_k''$ 
and $Y_k((x_j',x_j'',y_j)_{j \in \bN}) = y_k$).

Suppose that $\bP', \bP'' \in \Pc_\mu$.  
Write $\bP_z'(\cdot)$ (resp. $\bP_z''(\cdot)$) 
for the regular conditional probability of $\bP'\{X \in \cdot \, | \, Z=z\}$ 
(resp. $\bP''\{X \in \cdot \, | \, Z=z\}$).

Define a probability measure $\bQ$ on $\prod_{k \in \bN} (G_k \times G_k \times G_k)$
by
\[
\bQ\{(X',X'',Y) \in A' \times A'' \times B\}
=
\int_G \bP_z'(A') \bP_z''(A'') 1_B(z) \, (\bigotimes_{k \in \bN} \mu_k)(dz).
\]
By construction,
$\phi_k(X'_{k+1})^{-1} \, X'_k 
= \phi_k(X''_{k+1})^{-1} \, X''_k = Y_k$ for all $k \in \bN$,
$\bQ$-a.s., the distribution of the pair
$(X',Y)$ under $\bQ$ is the same as that of
the pair $(X,Z)$ under $\bP'$, and the distribution of 
the pair $(X'',Y)$ under $\bQ$ is the same as that of
the pair $(X,Z)$ under $\bP''$.  In particular, the distributions of
$X'$ and $X''$ under $\bQ$ are, respectively, $\bP'$ and $\bP''$.

Suppose for some $k \in \bN$ that 
$\Phi':G \to \bR$ and $\Phi'': G \to \bR$ are both bounded $\F_{k+1}^X$-measurable 
functions and $\Psi : G_k \to \bR$ is a bounded Borel function.
Then, 
$\Phi' \circ X' : \prod_{j \in \bN} (G_j \times G_j \times G_j)\to \bR$ 
is $\F_{k+1}^{X'}$-measurable and
$\Phi'' \circ X'' : \prod_{j \in \bN} (G_j \times G_j \times G_j)\to \bR$ 
is $\F_{k+1}^{X''}$-measurable, and hence, by the construction
of $\bQ$ (using the notations $\nu[\cdot]$ and
$\nu[\cdot \, | \, \cdot]$ for expectation and conditional expectation with respect
to a probability measure $\nu$),
\[
\begin{split}
\bQ[\Phi' \circ X' \, \Phi'' \circ X'' \, | \, \F^Y]
& =
\bQ[\Phi' \circ X'  \, | \, \F^Y] \; \bQ[\Phi'' \circ X'' \, | \, \F^Y] \\
& =
\bP'_Y[\Phi' \circ X] \; \bP''_Y[\Phi'' \circ X] \\
\end{split}
\]
is
$\F_{k+1}^Y$-measurable.  Thus, by the construction of $\bQ$
and the independence of the elements of the sequence 
$(Y_j)_{j \in \bN}$ under $\bQ$,
\[
\begin{split}
\bQ[\Phi' \circ X' \, \Phi'' \circ X'' \, \Psi \circ Y_k]
& =
\bQ[\bQ[\Phi' \circ X' \, \Phi'' \circ X'' \, \Psi \circ Y_k \, | \, \F^Y]] \\
& =
\bQ[\bQ[\Phi' \circ X' \, \Phi'' \circ X'' \, | \, \F^Y] \Psi \circ Y_k] \\
& =
\bQ[\bP'_Y[\Phi' \circ X] \,\bP''_Y[\Phi'' \circ X]] \; \bQ[\Psi \circ Y_k] \\
& =
\bQ[\Phi' \circ X' \, \Phi'' \circ X''] \; \bQ[\Psi \circ Y_k]. \\
\end{split}
\]
Therefore, by a standard monotone class
argument, $Y_k$ is independent of $\F_{k+1}^{(X',X'')}$.
Consequently, the sub-$\sigma$-fields
$\F_Y$ and $\F_\infty^{(X',X'')}$ are independent.

Suppose now that $\bP', \bP'' \in \Pc_\mu^{\mathrm{ex}}$.
Observe for $k < n$ that
\begin{equation}\label{coupling}
\begin{split}
& X_k' (X_k'')^{-1} \\
& \quad = 
\left[
\phi_k^n(X_n') \prod_{m=k}^{n-1} \phi_k^m(Y_m) \, Y_k
\right] 
\left[
\phi_k^n(X_n'') \prod_{m=k}^{n-1} \phi_k^m(Y_m) \, Y_k
\right]^{-1} \quad \bQ-\text{a.s.} \\
& \quad = \phi_k^n(X_n') \, \phi_k^n(X_n'')^{-1}, \\
\end{split}
\end{equation}
and so there exists a $G$-valued random variable $W \in \F^{X', X''}_{\infty}$ such that 
$W_k = X_k' (X_k'')^{-1}$, $\bQ$-a.s. From the above, $W$ is
independent of the sub-$\sigma$-field $\F_Y$.  By construction, $W$ takes values in the subgroup $H$.

Let $\bQ(\cdot \, | \, W = h)$ be the regular conditional probability for $\bQ$ 
given $W = h \in H$, so that
\begin{equation}
\label{convex_comb}
\bQ(\cdot) 
= \int_{H} \bQ(\cdot \, | \, W = h) \, \bQ\{W \in dh\}.
\end{equation}
It follows that
\[
\bQ\{X_k' = \phi_k(X_{k+1}') \, Y_k, \, \forall k \in \bN \, | \, W= h\} = 1
\]
for $\bQ\{W \in dh\}$-almost every $h \in H$. Moreover, because
$W$ is independent of $\F_Y$ it follows that 
$\bQ\{Y \in \cdot\} = \bQ\{Y \in \cdot \, | \, W = h\} = \bigotimes_{k \in \bN} \mu_k$
for $\bQ\{W \in dh\}$-almost every $h \in H$.
Thus, $\bQ\{X' \in \cdot \, | \, W = h\} \in \Pc_\mu$ for 
$\bQ\{\epsilon \in dh\}$-almost every $h \in H$ and, by \eqref{convex_comb},
\[
\bP'(\cdot) 
= \bQ\{X' \in \cdot\}
= \int_{H} \bQ\{X' \in \cdot \, | \, W = h\} \, \bQ\{W \in dh\}.
\]
This would contradict the extremality of $\bP'$
unless
\[
\bP'(\cdot) = \bQ\{X' \in \cdot \, | \, W = h\}, \;
\text{for 
$\bQ\{W \in dh\}$-almost every $h \in H$}.
\]
Similarly,
\[
\bP''(\cdot) = \bQ\{X'' \in \cdot \, | \, W = h\}, \;
\text{for 
$\bQ\{W \in dh\}$-almost every $h \in H$}.
\]
By \eqref{coupling}, 
\[
\bQ\{X_k' = h_k X_k'' \, \forall k \in \bN \, | \, W = h\} = 1, \;
\text{for 
$\bQ\{W \in dh\}$-almost every $h \in H$}.
\]
Therefore,
\[
\bP' = T_h(\bP''), \;
\text{for 
$\bQ\{W \in dh\}$-almost every $h \in H$}.
\]
\end{proof}

\begin{notation}
Given $\bP^0 \in \Pc_\mu^{\mathrm{ex}}$, let 
$H_\mu^{\mathrm{stab}}(\bP^0) := \{h \in H : T_h(\bP^0) = \bP^0\}$  be the stabilizer 
subgroup of the
point $\bP^0$ under the group action $(T_h)_{h \in H}$. 
\end{notation}

\begin{remark}
It follows from the transitivity of $H$ on $\Pc_\mu^{\mathrm{ex}}$
that for any two probability measures $\bP', \bP'' \in \Pc_\mu^{\mathrm{ex}}$
the subgroups $H_\mu^{\mathrm{stab}}(\bP')$ and $H_\mu^{\mathrm{stab}}(\bP'')$ 
are conjugate.
\end{remark}

\begin{corollary}
\label{nec_suff_unique}
A necessary and sufficient condition for $\# \Pc_\mu = 1$ is that
$H_\mu^{\mathrm{stab}}(\bP^0) = H$ for some, and hence all, 
$\bP^0 \in \Pc_\mu^{\mathrm{ex}}$.
\end{corollary}

\begin{proof}
This is immediate from Theorem~\ref{htransitive} 
and the observation that $\# \Pc_\mu = 1$
if and only if $\# \Pc_\mu^{\mathrm{ex}} = 1$.
\end{proof}

\begin{corollary}
If $H_\mu^{\mathrm{stab}}(\bP^0)$ is non-trivial for some, and hence all, 
$\bP^0 \in \Pc_\mu^{\mathrm{ex}}$, then $\Pc_\mu^{\mathrm{strong}} = \emptyset$.
\end{corollary}

\begin{proof} As we observed prior to the statement of Theorem~\ref{htransitive},
it is a consequence of that result that either
$\Pc_\mu^{\mathrm{strong}} = \emptyset$
or $\Pc_\mu^{\mathrm{strong}} = \Pc_\mu^{\mathrm{ex}}$.

Suppose that $\bP^0 \in \Pc_\mu^{\mathrm{strong}}$ is such that
$H_\mu^{\mathrm{stab}}(\bP^0)$ is non-trivial.  By working on an extended
probability space, we may assume that there is an
$H_\mu^{\mathrm{stab}}(\bP^0)$-valued random variable $(U_k)_{k \in \bN}$ 
that is independent
of $(X_k)_{k \in \bN}$ and is not almost surely constant. 
The distribution of the solution $(U_k X_k)_{k \in \bN}$ is also $\bP^0$ 
and, in particular, this solution is strong.  However, this implies that
\[
\begin{split}
\sigma(U_k X_k)
& \subseteq 
\sigma((\phi_j(U_{j+1}\, X_{j+1})^{-1} \, U_j \, X_j)_{j \ge k}) \\
& =
\sigma((\phi_j(X_{j+1})^{-1} \, X_j)_{j \ge k}) \\
& =
\F_k^Z \\
\end{split}
\]
for all $k \in \bN$, and hence $U_k$ is $\F_k^Z$-measurable
 for all $k \in \bN$,
because $X_k$ is $\F_k^Z$-measurable by the assumption that 
$\bP^0 \in \Pc_\mu^{\mathrm{strong}}$.  However, because the sequence
$(U_k)_{k \in \bN}$ is independent of the sequence
of $(X_k)_{k \in \bN}$ and not almost surely constant, it follows that
that $(U_k)_{k \in \bN}$ is not 
$\sigma((X_k)_{k \in \bN})$-measurable, and hence
{\em a fortiori}, $(U_k)_{k \in \bN}$ is not $\sigma((Z_k)_{k \in \bN})$-measurable.
We thus arrive at a contradiction.
\end{proof}

\section{Representation theory and the existence of strong solutions}
\label{S:representations_and_strong}

\begin{notation}
Let $\G$ be the set of all unitary, finite-dimensional representations of 
the compact group
$G = \prod_{k \in \bN} G_k$. 
\end{notation}

Any irreducible representations of $G$ is equivalent to a
tensor product representation of the form
\[
(g_k)_{k \in \bN} 
\mapsto 
\rho^{(k_1)}(g_{k_1}) \otimes \cdots \otimes \rho^{(k_n)}(g_{k_n}),
\]
where $\{k_1, \ldots, k_n\}$ is a finite subset of $\bN$
and $\rho^{(k_j)}$ is a (necessarily finite-dimensional)
irreducible representation of $G_{k_j}$
for $1 \le j \le n$.  Furthermore, an arbitrary element of $\G$ is
equivalent to a (finite) direct sum of irreducible representations.

\begin{notation}
For $k \in \bN$ write $\iota_k : G_k \mapsto G$ for the map that sends
$h \in G_k$ to $(e_1, \ldots, e_{k-1}, h, e_{k+1}, \ldots)$, where, as above,
$e_j$ is the identity element of $G_j$ for $j \in \bN$.
\end{notation}

Consider an arbitrary representation $\rho \in \G$. 
It is clear from the above that if $\bP \in \Pc_\mu^{\mathrm{strong}}$, 
then $\rho \circ \iota_k(X_k)$ is  $\F^Z_k$-measurable for all $k \in \bN$.
Note that $\rho \circ \iota_k$ is a representation of $G_k$
and all representations of $G_k$ arise this way.
On the other hand, because, by the Peter-Weyl theorem, the closure in the uniform norm
of the (complex) linear span of matrix entries of the irreducible representations
of $G_k$ is the vector space of continuous complex-valued functions on
$G_k$, it follows that if $\rho \circ \iota_k(X_k)$ is $\F^Z_k$-measurable for all $k \in \bN$
for an arbitrary representation $\rho \in \G$, then 
$\bP \in \Pc_\mu^{\mathrm{strong}}$.  This observation leads to the following
definition and theorem.

\begin{notation}
\label{Hc_strong}
Set
\[
\Hc_\mu^{\mathrm{strong}} 
:= 
\{ \rho \in \G:
\text{$\exists \bP \in \Pc_\mu^{\mathrm{ex}}$ 
such that 
$\rho \circ \iota_k (X_k)$ is $\F^Z_k$-measurable 
$\bP$-a.s. $\forall k \in \bN$}\}.
\]
\end{notation}

\begin{theorem}
\label{Hc_strong_nasc}
The set $\Pc_\mu^{\mathrm{strong}}$ of strong solutions is non-empty 
(and hence equal to $\Pc_\mu^{\mathrm{ex}}$) if and only if
$\Hc_\mu^{\mathrm{strong}} = \G$.
\end{theorem}

\begin{proof}
The result is immediate from the discussion preceding the statement of
the theorem once we note that if $\bP'$ and $\bP''$ both belong to
$\Pc_\mu^{\mathrm{ex}}$ then, by Theorem~\ref{htransitive}, 
there exists $h \in H$ such that $\bP''$ is the distribution
of $h X = (h_k X_k)_{k \in \bN}$ under $\bP'$ and so 
$\rho \circ \iota_k (X_k)$ is $\F^Z_k$-measurable $\bP''$-a.s. if and only if
$\rho \circ \iota_k(h_k X_k)$ is  $\F^Z_k$-measurable $\bP'$-a.s. (recall
that $Z_k = \phi(X_{k+1})^{-1} \, X_k = \phi(h_k X_{k+1})^{-1} \, h_k X_k$
when $h \in H$); therefore,
$\rho \circ \iota_k (X_k)$ is $\F^Z_k$-measurable $\bP''$-a.s. if and only
if $[\rho \circ \iota_k (h_k)] \, [\rho \circ \iota_k (X_k)]$
is $\F^Z_k$-measurable $\bP'$-a.s.,
which is in turn equivalent to $\rho \circ \iota_k (X_k)$ being  $\F^Z_k$-measurable $\bP'$-a.s.
by the invertibility of the matrix $\rho \circ \iota_k (h_k)$. Thus,
\[
\Hc_\mu^{\mathrm{strong}}
=
\{ \rho \in \G:
\text{
$\rho \circ \iota_k (X_k)$ is $\F^Z_k$-measurable $\bP$-a.s. $\forall k \in \bN$
}
\}
\]
for any $\bP \in \Pc_\mu^{\mathrm{ex}}$.
\end{proof}

Theorem~\ref{Hc_strong_nasc} is still somewhat unsatisfactory as
a criterion for the existence of strong solutions because it requires
a knowledge of the set $\Pc_\mu^{\mathrm{ex}}$ of extreme solutions.
We would prefer a criterion that was directly in terms of
the sequence $(\mu_k)_{k \in \bN}$.
In order to (partly) remedy this situation, we introduce the following objects.

\begin{notation}
\label{def_Rkl}
Fix $\rho \in \G$. For $k,\ell \in \bN$ with $k \le \ell$, set
\[
R_k^\ell 
:= 
\int_{G_\ell} \rho \circ \iota_k \circ \phi^{\ell}_{k}(z) \, \mu_\ell(dz).
\]
Let
\[
\Hc_\mu^{\mathrm{det}} 
:=
\{\rho \in \G :
\lim_{m \to \infty} \lim_{n \to \infty}
\left | \det (R_k^n R_k^{n-1} \cdots R_k^m) \right|
> 0 
\; \forall k \in \bN
\}
\]
and
\[
\Hc_\mu^{\mathrm{norm}} 
:= 
\{
\rho \in \G: 
\lim_{m \to \infty} \lim_{n \to \infty}
\| R_k^n R_k^{n-1} \cdots R_k^m\|
> 0 
\; \forall k \in \bN
\},
\]
where $\| \cdot \|$ is the $\ell^2$ operator norm on the appropriate space
of matrices.
\end{notation}

\begin{proposition}
\label{strongmeasurable} 
Fix $\bP \in \Pc_\mu$.
\begin{enumerate}
\item If $\rho \in \Hc_{\mu}^{\mathrm{det}}$, then
\[
\bP[\rho \circ \iota_k(X_k) \, | \, \F_{\infty}^{X} \vee \F_{k}^{Z}] 
= \rho \circ \iota_k(X_k) 
\]
for all $k \in \bN$.
In particular,
if $\bP \in \Pc_\mu^{\mathrm{ex}}$, 
then $\rho \circ \iota_k(X_k)$ is $\F_k^{Z}$-measurable 
for all $k \in \bN$.
\item If $\rho \notin \Hc_{\mu}^{\mathrm{norm}}$, then
\[
\bP[\rho \circ \iota_k(X_k) \, | \, \F_{\infty}^{X} \vee \F_{k}^{Z}] 
= 0
\]
for some $k \in \bN$.
In particular, 
if $\bP \in \Pc_\mu^{\mathrm{ex}}$, then $\rho \circ \iota_k(X_k)$ is not 
$\F_k^{Z}$-measurable for some $k \in \bN$.\\
\end{enumerate}
\end{proposition}

\begin{proof}
The proof follows that of an analogous result in
\cite{MR2365485} with modifications required
by the greater generality in which we are working.

Consider claim (i).  Fix $\rho \in \Hc_{\mu}^{\mathrm{det}}$ and $k \in \bN$.
For $\ell > k$ we have
\begin{equation}
\label{product_decomp}
\rho \circ \iota_k(X_k) 
= 
\rho \circ \iota_k \circ \phi_k^\ell(X_\ell) \,
\rho \circ \iota_k \circ \phi_k^{\ell-1}(Z_{\ell-1}) \cdots 
\rho \circ \iota_k \circ \phi_k^k(Z_k).
\end{equation}

For $k \le m \le n$ put
\[
\Xi_n^m 
:= 
\rho \circ \iota_k \circ \phi_k^n(Z_m) \cdots \rho \circ \iota_k \circ \phi_k^m(Z_m).
\]
Note that 
\[
\bP[\Xi_n^m] = R_k^n  \cdots R_k^m.
\]
For any $p \ge k$, the matrix $\rho \circ \iota_k \circ \phi_k^p$ is unitary,
and so  $\|\rho \circ \iota_k \circ \phi_k^p(h)\| = 1$ for all $h \in G_p$.
By Jensen's inequality, $\|R_k^p\| \le 1$. In particular, $|\det(R_k^p)| \le 1$.
Hence,
\[
\lim_{m \to \infty} \lim_{n \to \infty} |\det (\bP[\Xi_n^m]) |
\]
exists and is given by
\[
\sup_m \inf_{n \ge m} |\det(R_k^n)|  \cdots |\det(R_k^m)|.
\]
Moreover, there are constants $\epsilon > 0$ and $M \in \bN$
such that $| \det(\bP[\Xi_n^m]) | \ge \epsilon$ whenever
$n \ge m \ge M$.  It follows from Cramer's rule that the matrices
$\bP[\Xi_n^m]$ are invertible with uniformly bounded entries
for $n \ge m \ge M$.

Set $\Phi_n^m := \bP[\Xi_n^m]^{-1} \Xi_n^m$ for $n \ge m \ge M$.  
The matrices $\Phi_n^m$ have uniformly bounded entries
and
\[
\bP\left[\Phi_{n+1}^m \, | \, \sigma((Z_p)_{p=m}^n)\right] = \Phi_n^m,
\]
so that $(\Phi_n)_{n \ge m}$ is a bounded matrix-valued martingale with
respect to the filtration $(\sigma((Z_p)_{p=m}^n))_{n \ge m}$.  Thus,
$\lim_{n \to \infty} \Phi_n^m =: \Phi_\infty^m$ exists
and is $\F_m^Z$-measurable $\bP$-a.s. for each $m \ge M$.  Consequently, $\lim_{n \to \infty} \Xi_n^m =: \Xi_\infty^m$ also exists and is 
$\F_m^Z$-measurable $\bP$-a.s. for each $m \ge M$. Part (i) is now clear from \eqref{product_decomp}.

Now consider part (ii).
Fix $\rho \notin \Hc_\mu^{\mathrm{norm}}$ and $k \in \bN$ such that
\[
\lim_{m \to \infty} \lim_{n \to \infty}
\left\| R_k^n R_k^{n-1} \cdots R_k^m \right\|
=
0.
\]
It follows from \eqref{product_decomp} that for $n \ge m \ge k$
\[
\begin{split}
\bP\left[\rho \circ \iota_k(X_k) \, | \, \F_n^X \vee \sigma((Z_j)_{j=k}^m)\right] 
& = 
\rho \circ \iota_k \circ \phi_k^n(X_n)
R_k^{n-1} \cdots R_k^{m+1} \\
& \quad
\rho \circ \iota_k \circ \phi_m^k(Z_m) \cdots \rho \circ \iota_k \circ \phi_k^k(Z_k).
\\
\end{split}
\]
Since $\rho(g)$ is a unitary matrix for all $g \in G$, 
the norm of the right-hand side is at most $\|R_k^{n-1} \cdots R_k^{m+1}\|$,
which, by assumption, converges to $0$ as $n \to \infty$ followed by $m \to \infty$.
Thus, by the reverse martingale convergence theorem and the martingale
convergence theorem,
\[
\bP\left[\rho \circ \iota_k(X_k) \, | \, \F_{\infty}^{X} \vee \F_{k}^{Z}\right] 
=
\lim_{m \to \infty} \lim_{n \to \infty}
\bP\left[\rho \circ \iota_k(X_k) \, | \, \F_n^X \vee \sigma((Z_j)_{j=k}^m)\right] 
= 0.
\]
\end{proof}

The following result is immediate from Theorem~\ref{Hc_strong_nasc}
and Proposition~\ref{strongmeasurable}.

\begin{theorem}
\label{containment}
The following containments hold
\[
\Hc_\mu^{\mathrm{norm}} 
\supseteq \Hc_\mu^{\mathrm{strong}} 
\supseteq \Hc_\mu^{\mathrm{det}}.
\]
Thus,
$\Hc_\mu^{\mathrm{det}} = \G$ implies that $\Pc_\mu^{\mathrm{strong}} \ne \emptyset$ 
and
$\Hc_\mu^{\mathrm{norm}} \ne \G$ implies that $\Pc_\mu^{\mathrm{strong}} = \emptyset$.
\end{theorem}

The following is a straightforward equivalent of Theorem~\ref{containment} and
we omit the proof.

\begin{corollary}
\label{irreducible_restatement}
If 
\[
\lim_{m \to \infty} \lim_{n \to \infty}
\left| 
\det\left(
\prod_{\ell=m}^n \int_{G_\ell} \rho \circ \phi_k^\ell(z) \, \mu_\ell(dz)
\right) 
\right|
> 0
\]
for all irreducible representations $\rho$ of $G_k$ for all $k \in \bN$, then
$\Pc_\mu^{\mathrm{strong}} \ne \emptyset$.  If
\[
\lim_{m \to \infty} \lim_{n \to \infty}
\left\|\prod_{\ell=m}^n \int_{G_\ell} \rho \circ \phi_k^\ell(z) \, \mu_\ell(dz)) \right\|
= 0
\]
for some irreducible representation $\rho$ of $G_k$ for some $k \in \bN$, then
$\Pc_\mu^{\mathrm{strong}} = \emptyset$.
\end{corollary}

Under a further assumption, we get a representation theoretic
necessary and sufficient condition for the existence of strong solutions.

\begin{definition}
A Borel probability measure $\nu$ on a compact Hausdorff group $\Gamma$ is 
{\em conjugation invariant} if 
\[
\int_\Gamma f(g^{-1} x g) \, \nu(dx) = \int_\Gamma f(x) \, \nu(dx)
\]
for all $g \in \Gamma$ and bounded Borel functions $f:\Gamma \to \bR$.
\end{definition}

\begin{remark}
Note that if $\Gamma$ is abelian, then any Borel probability measure $\nu$
on $\Gamma$ is conjugation invariant.
\end{remark}

\begin{corollary}
\label{conjugation_invariant_equality}
Suppose that each probability measure $\mu_k$, $k \in \bN$, is conjugation
invariant. Then,
\[
\Hc_\mu^{\mathrm{norm}} 
= \Hc_\mu^{\mathrm{strong}} 
= \Hc_\mu^{\mathrm{det}}
\]
and
 $\Pc_\mu^{\mathrm{strong}} \ne \emptyset$
if and only if  each of these sets is $\G$ or, equivalently,
\[
\lim_{m \to \infty} \lim_{n \to \infty}
\left| \prod_{\ell=m}^n \int_{G_\ell} \chi \circ \phi_k^\ell(z) \, \mu_\ell(dz) \right|
> 0
\]
for each character $\chi$ of an irreducible representation of $G_k$ for all $k \in \bN$.
\end{corollary}

\begin{proof}
The result is immediate from Corollary~\ref{irreducible_restatement} and 
Lemma~\ref{characterssuffice} below.
\end{proof}

The following lemma is well-known, but we include a proof for 
the sake of completeness.

\begin{lemma}
\label{characterssuffice} 
If $\nu$ is a conjugation invariant Borel probability measure 
on a compact Hausdorff group $\Gamma$ and $\rho$ is an irreducible
representation of $\Gamma$ with character $\chi$, then
\[
\int_\Gamma \rho(x) \, \nu(dx) = \int_\Gamma \chi(x) \, \nu(dx) \times I,
\]
where $I$ is the identity matrix.
\end{lemma}

\begin{proof}
Let $\lambda$ be the normalized Haar measure on $\Gamma$.  By assumption,
\[
\int_\Gamma \rho(x) \, \nu(dx) 
= 
\int_\Gamma \int_\Gamma \rho(g^{-1} x g) \, \lambda(dg) \, \nu(dx).
\]
Now, for $x,y \in \Gamma$ we have
\[
\begin{split}
\int_\Gamma \rho(g^{-1} x g) \, \lambda(dg) \; \rho(y)
& =
\int_\Gamma \rho(g^{-1} x g y) \, \lambda(dg) \\
& =
\int_\Gamma \rho(y h^{-1} x h) \, \lambda(dh) \\
& =
\rho(y) \;  \int_\Gamma \rho(h^{-1} x h) \, \lambda(dh), \\
\end{split}
\]
and so the matrix $\int_\Gamma \rho(g^{-1} x g) \, \lambda(dg)$ commutes
with the matrix $\rho(y)$ for all $y \in \Gamma$.  It follows from Schur's Lemma
that $\int_\Gamma \rho(g^{-1} x g) \, \lambda(dg) = c I$ for some constant $c$,
and taking traces of both sides gives $c = \chi(x)$.
\end{proof}

\section{Freezing}
\label{S:freezing}

Recall that the Hilbert-Schmidt norm of a matrix $A$ is given by 
$\| A \|_{HS} := \tr(A^*A)^{\frac{1}{2}}$, where $A^*$
is the adjoint of $A$ (this norm is also called
the Frobenius norm and the Schur norm).  Write $d(\rho)$
for the dimension of a unitary representation $\rho \in \G$,
and note that $\| \rho(x)\|_{HS}^2 = \tr(I) = d(\rho)$.
If $\nu$ is a probability measure on $G$, then
$\| \int_G \rho(x) \, \nu(dx) \|_{HS}^2 \le d(\rho)$
by Jensen's inequality.

\begin{notation} Set
\[
\Hc_\mu^{\mathrm{freeze}} 
:= 
\left\{\rho \in \G :  
\sum_{m=k}^\infty
\left[ 
d(\rho) 
- \left \| \int_{G_k} \rho \circ \iota_k \circ \phi_k^m(z) \, \mu_m(dz) \right \|_{HS}^{2}
\right] 
< \infty \; \forall k \in \bN
\right\}.
\]
\end{notation}

\begin{proposition}
\label{freeze_equals_det}
The sets $\Hc_\mu^{\mathrm{freeze}}$ and $\Hc_\mu^{\mathrm{det}}$
are equal, and so $\Hc_\mu^{\mathrm{freeze}} = \Hc_\mu^{\mathrm{det}} = \G$ 
implies that $\Pc_\mu^{\mathrm{strong}} \ne \emptyset$. 
Moreover, if
each probability measure $\mu_k$, $k \in \bN$, is conjugation
invariant, then,
\[
\Hc_\mu^{\mathrm{norm}} 
= \Hc_\mu^{\mathrm{strong}} 
= \Hc_\mu^{\mathrm{det}}
= \Hc_\mu^{\mathrm{freeze}}
\]
and
 $\Pc_\mu^{\mathrm{strong}} \ne \emptyset$
if and only if each of these sets is $\G$ or, equivalently,
\[
\lim_{m \to \infty} \lim_{n \to \infty}
\left| \prod_{\ell=m}^n \int_{G_\ell} \chi \circ \phi_k^\ell(z) \, \mu_\ell(dz) \right|
> 0
\]
for each character $\chi$ of an irreducible representation of $G_k$ for all $k \in \bN$.
\end{proposition}

\begin{proof}
It suffices to show that $\Hc_\mu^{\mathrm{freeze}} = \Hc_\mu^{\mathrm{det}}$,
because the remainder of the result will then follow from
Theorem~\ref{containment}
and
Corollary~\ref{conjugation_invariant_equality}.

Fix $\rho \in \G$.  Write 
$0 \le \lambda_k^\ell(1) \le \cdots \le \lambda_k^\ell(d(\rho))$
for the eigenvalues of the matrix
\[
\left(\int_{G_k} \rho(z) \, \mu_k^\ell(dz)\right)^*
\left(\int_{G_k} \rho(z) \, \mu_k^\ell(dz)\right).
\]

Observe that
\[
\begin{split}
& \lim_{m \to \infty} \lim_{n \to \infty}
\prod_{\ell=m}^n \left| \det \int_{G_k} \rho(z) \, \mu_k^\ell(dz) \right| > 0\\
& \qquad \Longleftrightarrow \\
& \quad \lim_{m \to \infty} \lim_{n \to \infty}
\prod_{\ell=m}^n \left| \det \int_{G_k} \rho(z) \, \mu_k^\ell(dz) \right|^2 > 0 \\
& \qquad \Longleftrightarrow \\
& \quad \lim_{m \to \infty} \lim_{n \to \infty}
\prod_{\ell=m}^n \lambda_k^\ell(1) \cdots \lambda_k^\ell(d(\rho)) > 0 \\
& \qquad \Longleftrightarrow \\
& \quad \sum_{m=k}^\infty \left[(1 - \lambda_k^m(1)) + \cdots + (1 - \lambda_k^m(d(\rho)))\right] < \infty \\
& \qquad \Longleftrightarrow \\
& \quad \sum_{m=k}^\infty \left[d(\rho) - \left( \lambda_k^m(1)) + \cdots + \lambda_k^m(d(\rho))\right)\right] < \infty \\
& \qquad \Longleftrightarrow \\
& \quad \sum_{m=k}^\infty
\left[ 
d(\rho) 
- \left \| \int_{G_k} \rho \circ \iota_k \circ \phi_k^m(z) \, \mu_m(dz) \right \|_{HS}^{2}
\right] 
< \infty, \\
\end{split}
\]
as required.
\end{proof}

Given Proposition~\ref{freeze_equals_det}, the reader may wonder why
we introduced the set $\Hc_\mu^{\mathrm{freeze}}$.  The equivalence
established in Proposition~\ref{freeze_equals_det} makes the
proof of the following result considerably more transparent.

\begin{proposition} 
\label{freeze_meaning}
Suppose that each group $G_k$, $k \in \bN$, is finite.
Then, 
$\Hc_\mu^{\mathrm{det}} = \Hc_\mu^{\mathrm{freeze}} = \G$ if and only if
for some (equivalently, all) $\bP \in \Pc_\mu$
there are constants $c_{k,m} \in G_k$, $k,m \in \bN$,
$k \le m$, such that 
\[
\bP\{\phi_k^m(Z_m) \ne c_{k,m} \, \mathrm{i.o.}\} = 0
\]
for all $k \in \bN$.
\end{proposition}

\begin{proof}
Write $\mu_k^m$ for the probability measure on $G_k$ that is the push-forward
of the probability measure $\mu_m$ on $G_m$ by the map $\phi_k^m: G_m \to G_k$.
For simplicity, we write $\mu_k^m(g)$ instead of $\mu_k^m(\{g\})$ for $g \in G_k$.
It is clear that
$\bP\{\phi_k^m(Z_m) \ne c_{k,m} \, \mathrm{i.o.}\} = 0$
$k \le m$ for all $k \in \bN$
for some family of constants $c_{k,m} \in G_k$, $k,m \in \bN$,
if and only if
$\bP\{\phi_k^m(Z_m) \ne c_{k,m}^* \, \mathrm{i.o.}\} = 0$
where $c_{k,m}^*$ is any family with the property
\[
\mu(c_{k,m}^*) = \max\{\mu_k^m(g) : g \in G_k\}
\]
and, by the Borel-Cantelli lemma, this in turn occurs if and only if
\[
\sum_{m=k}^\infty \mu(G_k \backslash \{c_{k,m}^*\}) < \infty
\]
for all $k \in \bN$.

Now,
\[
\left( \sum_{g \in G_k} \mu_k^m(g)^2 \right)^{1/2} 
\geq 
\max_{g \in G_k} \mu_k^m(g) 
=
\mu_k^m(c_{k,m})
= 
\mu_k^m(c_{k,m}) \sum_{g \in G_k} \mu_k^m(g) 
\geq 
\sum_{g \in G_k} \mu_k^m(g)^2.
\]
By Parseval's equality,
\[
\sum_{g \in G_k} \mu_k^m(g)^2 
= 
\frac{1}{\# G_k} \sum_{\rho \in \hat{G_k}} 
d(\rho) \left \| \sum_{g \in G_k} \rho(g) \mu_k^m(g) \right \|_{HS}^{2},
\]
and hence
\[
\begin{split}
& 1 - 
\left( 
\frac{1}{\# G_k} \sum_{\rho \in \hat G_k} d(\rho)\left\| \sum_{g \in G_k} \rho(g) \mu_k^m(g) \right\|_{HS}^{2} \right) \\
& \quad \geq
\mu_k^m(G_k \backslash\{c_{k,m}\}) \\
& \quad \geq 
1 - 
\left( 
\frac{1}{\# G_k} \sum_{\rho \in \hat G_k} d(\rho) \left\| \sum_{g \in G_k} \rho(g) \mu_k^m(g) \right\|_{HS}^{2} \right)^{1/2}. \\
\end{split}
\]

Note for a sequence of constant $(a_n)_{n \in \bN} \subset [0,1]$ that 
$\sum_{n \in \bN} (1 - a_n) < \infty$ if and only if 
$\sum_{n \in \bN} (1 - a_n^2) < \infty$.  Note also that
\[
1 = \frac{1}{\# G_k} \sum_{\rho \in \hat G_k} d(\rho)^2.
\]
Thus,
\[
\sum_{m=k}^\infty \mu(G_k \backslash \{c_{k,m}^*\}) < \infty
\]
for all $k \in \bN$ if and only if 
\[
\sum_{m=k}^\infty 
\frac{1}{\# G_k} \sum_{\rho \in \hat G_k} 
d(\rho) \left[d(\rho) - \left\| \sum_{g \in G_k} \rho(g) \mu_k^m(g) \right\|_{HS}^{2}\right]
< \infty
\]
for all $k \in \bN$,
which is in turn equivalent to 
\[
\sum_{m=k}^\infty 
\sum_{\rho \in \hat G_k} 
\left[d(\rho) - \left\| \sum_{g \in G_k} \rho(g) \mu_k^m(g) \right\|_{HS}^{2}\right]
< \infty
\]
for all $\rho \in \hat G_k$ for all $k \in \bN$.

A decomposition of the representation $\rho \circ \iota_k$ of $G_k$
for some $\rho \in \G$ into irreducibles shows that the last
condition is equivalent to the one in the statement.
\end{proof}

\begin{remark}
It follows from Proposition~\ref{freeze_equals_det} 
and Proposition~\ref{freeze_meaning} that if
each group $G_k$, $k \in \bN$, is finite 
and for some (equivalently, all) $\bP \in \Pc_\mu$
there are constants $c_{k,m} \in G_k$, $k,m \in \bN$,
$k \le m$, such that 
\[
\bP\{\phi_k^m(Z_m) \ne c_{k,m} \, \mathrm{i.o.}\} = 0
\]
for all $k \in \bN$, then $\Pc_\mu^{\mathrm{strong}} \ne \emptyset$.
Moreover, these two conditions are equivalent when
each probability measure $\mu_k$, $k \in \bN$, is conjugation
invariant.  Also, for the special case when $G_k = \Gamma$,
$k \in \bN$, for some fixed finite group $\Gamma$ and each homomorphism
$\phi_k: \Gamma \to \Gamma$ is the identity, it follows from
Corollary 2.6 of \cite{MR2653259} that the two conditions are equivalent.
It would be interesting to know the status of the reverse 
implication in general.
\end{remark}

\section{Groups indexed by the lattice}
\label{S:lattice}

Recall from the Introduction the example of our general set-up where
$G_k := G_{1,k} \times G_{2,k-1} \cdots \times G_{k,1}$ with
each group $G_{i,j}$ a copy of some fixed compact abelian
group $\Gamma$ and the homomorphism
$\phi_k$ is given by
\[
\phi_k(g_{1,k+1}, g_{2,k}, \ldots, g_{k+1,1}) 
:= (g_{1,k+1} + g_{2,k}, g_{2,k}+g_{3,k-1}, \ldots, g_{k,2} + g_{k+1,1}).
\]
We will consider the particular case where $\Gamma$ is
$\mathbb{Z}_{p}$, the group of integers modulo some prime number $p$.

Because $\mathbb{Z}_p$ is abelian, all its irreducible representations of $G$
are one-dimensional.  The irreducible representations are the trivial one
and those of the form $\rho(g) = \prod_{n=1}^{m} \exp\left(\frac{2\pi i z_n}{p} g_{i_n, j_n}\right)$ for
some $m$, pairs $(i_1, j_1), \ldots, (i_m, j_m) \in \bN^2$, and $1 \leq z_n \leq p-1$.

The homomorphism $\phi^\ell_k$ maps $(g_{1, \ell}, \ldots, g_{\ell, 1}) \in
G_{\ell}$ to $(h_{1,k}, \ldots, h_{k,1}) \in G_k$ where
\[
h_{i,k+1-i} = \sum_{j=0}^{\ell-k} \binom{\ell-k}{j} g_{i+j,\ell+1-i-j} \in \mathbb{Z}_p.
\]

Set $f(m,n) := \binom{m}{n} \mod{p}$. When we restrict to $G_k$, the representation 
$\rho \circ \iota_k$ is of the form 
$\prod_{i=1}^{k} \exp\left(\frac{2\pi z_i}{p} g_{i, k+1-i}\right)$ 
with $0 \le z_i \le p-1$. 
We therefore need to evaluate
\[
R_k^\ell
=
\int_{G_{\ell}} 
\prod_{i=1}^{k} \prod_{j=0}^{\ell-k} 
\exp\left(\frac{2\pi z_i}{p} f(\ell-k,j) g_{i+j, \ell+1-i-j}\right) 
\, \mu_{\ell}(dg_{\ell})
\]
to determine whether or not $\Pc_\mu^{\mathrm{strong}} = \emptyset$.
The following theorem of Lucas (see \cite{MR1483922}) gives the value of $f$.

\begin{theorem} 
\label{Lucas}
Let $m, n$ be non-negative integers and $p$ a prime number. Suppose
\[m = m_k p^k + \ldots + m_1 p + m_0\]
and
\[n = n_k p^k + \ldots + n_1 p + n_0.\]
Then,
\[
\binom{m}{n} 
= 
\prod_{i=0}^{k} \binom{m_i}{n_i} 
\mod{p}. 
\] 
Equivalently, if $m_0$ and $n_0$ are the least non-negative residues of 
$m$ and $n$ mod $p$, then 
$\binom{m}{n} = \binom{\lfloor{m/p}\rfloor}{\lfloor{n/p}\rfloor}\binom{m_0}{n_0}$.
\end{theorem}

Rather than use Theorem~\ref{Lucas} directly to construct interesting
examples, we consider a consequence of it for the case $p=2$.
Suppose that
$\mu_k = \mu_{1,k} \otimes \cdots \otimes \mu_{k,1}$ where
$\mu_{i,k+1-i}\{1\} = \pi_k = 1 - \mu_{i,k+1-i}\{0\}$ for
some $0 \le \pi_k \le 1$.  

Define $x = (x_{m, \ell + 1 - m})_{m=1}^{\ell} \in G_{\ell}
= G_{1,\ell} \times \cdots \times G_{\ell ,1} \cong \mathbb{Z}_2^{\ell}$
by 
\[
x := \sum_{i=1}^{k} \sum_{j=0}^{\ell-k} z_i f(\ell-k,j) e^{(i+j, \ell+1-i-j)},
\]
where the arithmetic is performed modulo $2$ and 
$e^{(m, \ell+1-m)} \in G_{\ell}$ is the vector with
$e_{m, \ell+1-m}^{(m, \ell+1-m)} = 1$ and 
$e_{n, \ell+1-n}^{(m, \ell+1-m)} = 0$ for $n \ne m$.
Then,
\[
\int_{G_{\ell}} 
\prod_{i=1}^{k} \prod_{j=0}^{\ell-k}
\exp\left(\frac{2\pi z_i}{p} f(\ell-k,j) g_{i+j, \ell+1-i-j}\right) 
\, \mu_{\ell}(dg_{\ell})
=
(1 - 2 \pi_{\ell})^{M(k,\ell,z)},
\]
where 
\[
M(k,\ell,z)
:=
\#\{1 \le m \le \ell : x_{m, \ell+1-m} = 1\}.
\]

Observe that if $x_{m, \ell+1-m} = 1$, then 
\[
\sum_{j=0}^{\ell-k} f(\ell-k,j) e_{m, \ell+1-m}^{(i+j, \ell+1-i-j)} = 1
\]
for some $1 \le i \le k$ with $z_i = 1$.  Now
\[
\begin{split}
& \#
\{
1 \le m \le \ell 
: 
\sum_{j=0}^{\ell-k} f(\ell-k,j) e_{m, \ell+1-m}^{(i+j, \ell+1-i-j)} = 1
\} \\
& \quad =
\#
\{
1 \le m \le \ell 
: 
f(\ell-k, m - i)  = 1, \, i \le m \le i+\ell-k
\} \\
& \quad =
\#
\{
i \le m \le i + \ell-k
: 
f(\ell-k, m - i)  = 1
\} \\
& \quad =
\#
\{
0 \le m \le \ell-k
: 
f(\ell-k, m)  = 1
\}. \\
\end{split}
\]

As remarked in \cite{MR1483922}, a consequence of the following theorem
of Kummer from 1852 that the number of the binomial coefficients $\binom{m}{n}$,
$0 \le n \le m$, which are odd is $2^{N(m)}$, where $N(m)$ is the number
of times that the digit $1$ appears in the base $2$ representation of $m$.

\begin{theorem}
Let $m, n$ be non-negative integers and $p$ a prime number.
The greatest power of $p$ that divides
$\binom{m}{n}$ is given by the number of ``carries'' that are necessary when we 
add $m$ and $n-m$ in base $p$.
\end{theorem}

Thus,
\[
M(k,\ell,z) 
\le 
k 2^{N(\ell-k)}
\]
and $M(k,\ell,z) = 2^{N(\ell-k)}$ when $\#\{1 \le i \le k : z_i = 1\} = 1$.

Therefore, if we assume $\pi_n \to 0$ as $n \to \infty$,
then we are interested in whether
\[
\lim_{\ell \to \infty}
\prod_{r=1}^\ell (1 - 2 \pi_{h+r})^{2^{N(r)}} \ne 0
\]
for all $h \in \bN$ or, equivalently, whether
\[
\sum_{r=1}^\infty 2^{N(r)} \pi_{h+r} < \infty
\]
for all $h \in \bN$.

For example, fix a positive integer $a$ and an increasing function 
$b: \bN \to \bN$ 
such that $a \le b(m) < m$ and $\lim_{m \to \infty} b(m) = \infty$.
Suppose that $\pi_n = 0$ unless $2^m + 2^{b(m)} - 2^a \le n \le 2^m + 2^{b(m)}$ 
for some $m \in \bN$.  Note for any $h \in \bN$ 
that
\[
\sum_{r=1}^\infty 2^{N(r)} \pi_{h+r}
=
\sum_{s=k+1}^\infty 2^{N(s-h)} \pi_s
\]
and this sum is finite if and only if
\[
\sum_{n=1}^\infty 2^{b(\log_2 n)} \pi_n
\]
is finite.

Thus,
$\Pc_\mu^{\mathrm{strong}} \ne \emptyset$ if and only if
$\sum_{n=1}^\infty 2^{b(\log_2 n)} \pi_n < \infty$
in this case.  On the other hand,
$\bP\{Z_k \ne 0 \, \mathrm{i.o.}\} > 0$ (equivalently, 
$\bP\{Z_k \ne 0 \, \mathrm{i.o.}\} = 1$) if and only if
$\sum_{n=1}^\infty n \pi_n < \infty$.
Therefore, when $\lim_{m \to \infty} m - b(m) = \infty$ it is possible
to construct $(\pi_n)_{n \in \bN}$ such that almost surely
infinitely many ``bits'' are ``corrupted'' 
and yet strong solutions still exist.

\section{Automorphisms of the Torus}
\label{S:torus}

Consider the torus group
$\bT^2 = \mathbb{R}^{2} / \mathbb{Z}^{2}$. We  write an
element $x \in \bT^2$ as a column vector $x = (x_1, x_2)^\top \in [0,1)^{2}$, 
where $\top$ denotes the transpose of a vector.

Any $2 \times 2$ $\bZ$-valued matrix $S$ defines a homomorphism
$x \mapsto S x$ from $\bT^2$ to itself if we do ordinary matrix
multiplication modulo $\bZ^2$.  If the matrix $S$ has determinant
$1$, then this homomorphism is invertible.  Such a transformation
is called a \emph{linear toral automorphism}.
 
Note that if
\[
S =
\begin{pmatrix} a & b \\ c & d
\end{pmatrix},
\]
then the eigenvalues of $S$ are 
\[
\frac{1}{2} (a + d \pm \sqrt{a^2 + 4 b c - 2 a d + d^2})
=
\frac{1}{2} (a + d \pm \sqrt{(a+d)^2 - 4}),
\]
Thus, the eigenvalues are real and distinct unless
$a+d$ is $0$, $\pm 1$ or $\pm 2$, in which case 
the pairs of eigenvalues are, respectively
$\{\pm i\}$, 
$\{\frac{1}{2}(1 \pm i \sqrt{3})\}$,
$\{\frac{1}{2}(-1 \pm i \sqrt{3})\}$, 
$\{1,1\}$,
and $\{-1,-1\}$.
Note that in each of the latter cases
the eigenvalues lie on the unit circle.

\begin{definition} 
A \emph{ergodic toral automorphism} is a linear toral automorphism given by a 
matrix $S$ with no eigenvalues on the unit circle.
\end{definition}

For some of the
more probabilistic properties of ergodic toral automorphisms, see \cite{MR0294602}. 
Such mappings are the prototypical examples of Anosov systems that
have been the subject of intensive study dynamical systems world (see \cite{MR0253352}).

A hyperbolic linear toral automorphism has two real eigenvalues 
$\lambda_1 > 1 > \lambda_1^{-1} = \lambda_2$.
These eigenvalues are irrational and the corresponding 
(right) eigenvectors $v^1$ and $v^2$ 
have irrational slope (see, for example Section 5.6 of
\cite{MR1249274}).

\begin{theorem} Suppose for every $i \in \bN$ that the group $G_i$ is a copy of $\bT^2$
and that the homomorphism $\phi_i$ is a fixed 
ergodic toral automorphism given by a matrix $S$. Suppose the noise 
distribution $\mu_k$ is a fixed measure 
$\mu^*$ that
satisfies $\mu^*(A) \geq \epsilon \lambda (A \cap B)$ for every Borel set
$A$, where $\epsilon > 0$, $\lambda$ is normalized Haar measure,
and $B$ is a fixed Borel set $B$ with $\lambda(B)>0$. 
Then, $\Pc_\mu^{strong} = \emptyset$. 
\end{theorem}

\begin{proof}

We need to evaluate $R^{\ell}_k = \int_{\bT^2} \rho \cdot \iota_k \cdot \phi_k^{\ell}(z) \mu_{\ell}(dz)$. Let $\nu$ be the measure defined
by $\nu(A) = \epsilon \lambda(A \cap B)$ a Borel set $A$, where $\epsilon$,
$\lambda$ and $B$ are as in the statement.
Observe that
\begin{align*}
|R^{\ell}_k| &\leq \int_{\bT^2G_\ell} |\rho \cdot \iota_k \cdot \phi_k^{\ell}(z)| \, (\mu_{\ell} - \nu)(dz) + \int_{\bT^2} |\rho \cdot \iota_k \cdot \phi_k^{\ell}(z)| \, \nu(dz)| \\
&\leq \int_{\bT^2} \, (\mu_{\ell} - \nu)(dz) + \left| \int_{\bT^2} \rho \cdot \iota_k \cdot \phi_k^{\ell}(z) \, \nu(dz) \right|,
\end{align*}
and note that the last term on the right-hand side is 
$\left| \int_{\bT^2} \rho \cdot \iota_k(z) \, (\nu \cdot \phi_{k}^{\ell})^{-1}) (dz) \right|$.

As noted in Section 5.6 of \cite{MR1249274}, any ergodic toral automorphism $S$ 
exhibits \emph{topological mixing}: for any Borel sets $A, B \subseteq \mathbb{R}^2$, 
$\lim_{n \to \infty} \frac{\lambda(S^n B) \cap A}{\lambda(B)} = \lambda(A)$. 
Because $\phi_{k}^{\ell}$ is a ergodic toral automorphism, so is $(\phi_{k}^{\ell})^{-1}$. 
Therefore, 
$\lim_{\ell \to \infty} \left| \int_{\bT^2} \rho \cdot \iota_k(z) (\nu \cdot \phi_{k}^{\ell})^{-1} (dz) \right| 
= \left| \int_{\bT^2} \rho \cdot \iota_k(z) \epsilon \lambda (dz) \right| = 0$. 
Consequently, $|R^{\ell}_k| \leq \int_{\bT^2} (\mu_{\ell} - \nu)(dz) = 1 - \epsilon \lambda(B)$ 
for every non-trivial representation $\rho$, and hence
\[
\lim_{m \to \infty} \lim_{n \to \infty}
| R_k^n R_k^{n-1} \cdots R_k^m|
= 0 
\; \forall k \in \bN,
\]
showing that $\Pc_\mu^{strong} = \emptyset$.
\end{proof}

Every finite-dimensional unitary representation of $G_i$ is of the form,
\[
x \mapsto e^{2 \pi i(z \cdot x)},
\]
where $z$ is a vector $(z_1, z_2) \in \bZ^{2}$
and $z \cdot x$ is the usual inner product.
Hence, if we lift this representation to a representation of $G$ we have
\[
R_k^\ell = \int_{\bT^2} e^{2 \pi i(z \cdot S^{\ell-k} x)} \, \mu_\ell(dx).
\]

Suppose that the probability measure 
$\mu_\ell$ is concentrated on the set of multiples of
the eigenvector $v^2$ associated with the eigenvalue
$\lambda_2 \in (0,1)$.  Then,
\[
R_k^\ell 
= 
\int_\bR e^{2 \pi i(t \lambda_2^{\ell-k} z \cdot  v^2)} \, \nu_\ell(dt) 
\] 
for some probability measure $\nu_\ell$ on $\bR$.  It is clear that
under appropriate hypotheses 
\[
\lim_{m \to \infty} \lim_{n \to \infty}
| R_k^n R_k^{n-1} \cdots R_k^m|
> 0 
\; \forall k \in \bN
\]
and hence, by
Corollary~\ref{irreducible_restatement},
$\Pc_\mu^{\mathrm{strong}} \ne \emptyset$.
For example, if $\nu_\ell = \nu$ for all $\ell \in \bN$
for some fixed probability measure $\nu$ on $\bR$, then it
suffices that $\int_\bR |t| \, \nu(dt) < \infty$.  In particular,
it is possible to construct examples where $\mu_1 = \mu_2 = \ldots$ 
is a measure that has all of $\bT^2$ as its closed support and yet
$\Pc_\mu^{\mathrm{strong}} \ne \emptyset$.

\def\cprime{$'$}
\providecommand{\bysame}{\leavevmode\hbox to3em{\hrulefill}\thinspace}
\providecommand{\MR}{\relax\ifhmode\unskip\space\fi MR }
\providecommand{\MRhref}[2]{%
  \href{http://www.ams.org/mathscinet-getitem?mr=#1}{#2}
}
\providecommand{\href}[2]{#2}

\end{document}